\documentclass[a4paper]{article} %
\usepackage{graphicx,amssymb} %
\usepackage{amsmath,amsthm}
\usepackage{amssymb}
\usepackage[T1]{fontenc}
\usepackage{amsmath}
\usepackage{amsfonts}
\usepackage{amssymb}
\usepackage{amsthm}
\usepackage{graphicx}
\usepackage{makeidx}
\usepackage{color}
\newcommand{\C}{\mathbb{C}}
\newcommand{\R}{\mathbb{R}}
\newcommand{\N}{\mathbb{N}}
\newcommand*{\factorial}[1]{ #1 !}
\newtheorem{thm}{Theorem}[section]
\newtheorem{cor}[thm]{Corollary}
\newtheorem{lem}[thm]{Lemma}
\newtheorem{prop}[thm]{Proposition}
\newtheorem{prob}[thm]{Problem}
\newtheorem{defn}[thm]{Definition}

\theoremstyle{definition}

\newtheorem{rem}[thm]{Remark}
\newtheorem{exa}[thm]{Example}
\newtheorem{pre}[thm]{Proof}

\numberwithin{equation}{section}

\usepackage{enumerate}
\begin{document}


\baselineskip=17pt


\title{On some quasi-analytic classes }

\author{ \bf{ Abdelhafed Elkhadiri}\\
University Ibn Tofail,
Faculty of Sciences\\
Kenitra, Morocco\\
E-mail: elkhadiri.abdelhafed@uit.ac.ma}
\date{}

\maketitle


\renewcommand{\thefootnote}{}



\renewcommand{\thefootnote}{\arabic{footnote}}


\begin{abstract}
	Using the so called monotonicity property, we prove that 
	The Borel mapping restricted to
some quasi-anlytic classes	is never onto.

\end{abstract}

\section{Introduction}
Analytic functions, on an  interval $ [a,b]\subset \R$, $ a < b$, possess the following two equivalent properties:
\begin{enumerate}
\item [$\mathcal{B}$)] An analytic  function  is determined on $ [a,b]$ as soon as it is known in a subinterval of $ [a,b]$.
\item [$\mathcal{DC}$)] An analytic  function is determined  on $ [a,b]$  by its value and the values of its successive derivatives at a point $ c\in [a,b]$.
\end{enumerate}
This is a direct consequence of the fact that analytic functions  are developable in a Taylor's series in the neighborhood of the point and the Taylor's development identifies the function.
It was thought for a long time that the analytic functions where the only ones which were determined by their values and the values of there derivatives at a single point of the interval $ [a,b]$. It was Borel who first proved the existence of functions belonging to more general classes than that of analytic functions, which possessed the property that they are determined by their values and the values of all their successive derivatives at a single point $x_0\in [a,b]$, see
\cite{Borel}. Hi gave them the name {\it{quasi-analytic}} functions.\\
For any subvector space $\mathcal{A}$ of the ring of $ C^\infty$ functions on the interval $[a,b]$, $C^\infty( [a,b])$,
 if condition  $\mathcal{DC})$ is verified by the elements of $\mathcal{A}$,  then  condition $\mathcal{B})$ is  verified. In other words in the sitting of $C^\infty$ functions,  condition $\mathcal{DC})$ is stronger than condition $\mathcal{B})$.  The reciprocal is not true in general in the class of $ C^\infty$ functions. Take the function defined in the interval $[0,1]$ by $ f(0)= 0$ and $f(x)= e^{-\frac{1}{x}}$ if $x\in ]0,1]$.
 \\ The condition $\mathcal{B})$  is the definition of quasi-analyticity that Bernstein had adopted, see \cite{Bern}, while condition $\mathcal{DC})$ was adopted by Denjoy-Carleman.
Bernstein's definition of quasi-analyticity does not use the fact that the functions are differentiable, and can therefore be defined even for continuous functions.\\ In this paper we will say that the sub vector space  $\mathcal{A}\subset C^\infty( [a,b])$ is quasi-analytic, if $\mathcal{A}$ satisfied condition $\mathcal{DC})$.

Let's mention a second property that will interest us in this paper, which we call {\it{analytic monotonicity property}}.\\ If $f$ is an analytic functions on the interval $ [a,b]\subset \R$ and $ c\in [a,b[$, suppose that $ f^{(n)}(c) \geq 0$, $\forall n\in\N$, then there exists $ \eta >0$ such that $ [c,\eta]\subset [a,b[$   and $$\forall x\in [c,\eta], \forall n\in\N, \,\,f^{(n)}(x) \geq 0.$$ We are interested to see if the monotonicity property remains valid for these classes of quasi-analytic functions.
The importance of  the monotonicity property for these quasi-analytic classes lies in the following fact:\\
 We denote by $\R[[x]]$ the ring of  formal series in one variable with real coefficients. If $ c\in [a,b]$, the  Borel mapping $$T_c: C^\infty([a,b])\to \R[[x]]$$  is the function that associates to each function $ f\in C^\infty([a,b])$, the   series $\sum\limits_{n=0}^\infty \frac{f^{(n)}(c)}{n!}x^n$. By definition, the restriction of the Borel mapping to any quasi-analytic class is injective. It is a classical result due to Carleman, \cite{Car},  \cite{Car1}, that the restriction of the Borel mapping to the quasi-analytic class of Denjoy-Carleman, see below, is never onto.
We will see that this theorem remains true for other quasi-analytic classes with the property of monotonicity.
Many authors have investigated Carlman's proof by using techniques from functional analysis, see \cite{Thie}, theorem 3.
The  proof giving here is  direct and does not use functional analysis techniques.

\section{ Quasi-analytic functions of a real variable.}
In this section we are concerned with the approach used by Denjoy  to construct some classes of quasi-analytic real functions.
\subsection{Quasi-analytic  functions according to Denjoy's point of view.}
 Let $f$  be a $C^\infty$ function on the interval $[a,b]$ and put, for each $n\in\N$, $$ M_n(f) = \sup\limits_{x\in [a,b]}|f^{(n)}(x)|,$$
 where $f^{(n)}$ is the nth derivative of the function $f$.

Recall the characterization, given by Pringsheim, of analytic functions on the interval $[a,b]$. If
$f$ is a $C^\infty$  function on $[a,b]$, then
 \begin{equation}
f\,\,\,\mbox{ is analytic on the interval}\,\,\,\,[a,b]     \Leftrightarrow\,\,\sqrt[n]{M_n(f)}\leq C. n,\,\,\,\,\,\,\,\forall n\in\N,
\end{equation}
where $C$ is a positive constant independent of the integer $n$.\\
Denjoy asked himself  if it was not possible to enlarge the class of analytic functions on $[a,b]$ without losing the condition $\mathcal{DC})$. For this purpose he defined different classes of $C^\infty$ functions on the interval $[a,b]$ characterized by conditions

 \begin{equation}
\sqrt[n]{M_n(f)}\leq C\,n\,\ln n,\,\,\,\,\,\forall n > 1
\end{equation}

\begin{equation}
\sqrt[n]{M_n(f)}\leq C\,n\,\ln(n)\,\ln(\ln n),\,\,\,\,\,\forall n > e
\end{equation}
and so on.\\
He proved  that the functions of this different classes were still satisfied the condition $\mathcal{DC})$. Denjoy noticed that the reciprocal of the second members of the inequalities (2.2) and (2.3) are the general terms of divergent series. He was therefore led to announce the following theorem without prove it.
\begin{thm}
Let $f$  be a $C^\infty$ function on the interval $[a,b]$. The function is completely determined in the whole interval $[a,b]$ by its value and the values of its derivatives in any point of $[a,b]$, if the series of positive terms:
\begin{equation}
 \frac{1}{M_1(f)} + \frac{1}{\sqrt[2]{M_2(f)}} + \frac{1}{\sqrt[3]{M_3(f)}}+ \ldots + \frac{1}{\sqrt[n]{M_n(f)}}+\ldots
 \end{equation}
 is divergent
\end{thm}

In order to show this result, Carleman considered classes of functions more general than those considered by Denjoy. He proceeds as follows:\\
Let $M=(M_n)_{n\in\N}$ be a sequence of positive numbers. We denote by $ C_M([a,b])\subset C^\infty([a,b])$ the class of infinitely differentiable functions on the interval $[a,b]$ satisfying
$$ |f^{(n)}(x)|\leq c. C^n M_n,\,\,\,\,\forall  n\in \N,\,\,\,\forall x\in [a,b],$$
 where $c,\,\,C$ are positive constants (depending on $f$, but not on $n$). \\
We remark that if  $c$ is omitted in the definition, then if $n=0$ we have $\sup\limits_{x\in [a,b)}|f(x)|\leq  M_0$ which is restrictive. In the following, we suppose $ M_0 = 1$. The class $ C_M([a,b])$ is a vector space.\\
A class that satisfies condition $\mathcal{DC}$) will be called quasi-analytic class in the sense of Denjoy-Carleman.
\begin{prop}
Let $M=(M_n)_{n\in\N}$ be a sequence of positive numbers.  The class $ C_M([a,b])$ is quasi-analytic if  and only if condition $\mathcal{B}) $ is satisfied.
\end{prop}
\begin{pre}
It is enough to show that condition $\mathcal{B}) $ implies that the class is quasi-analytic.
Let $ f\in C_M([a,b])$ such that $ f^{(n)}(c)= 0,\,\,\,\forall n\in\N$, where $c\in ]a,b[$. Let g be the function defined as follows:
$g (x) = 0$  if $x\in [a, c]$ and $g (x) = f(x)$ if $x\in [c, b]$. \\
It is a straight  forward observation that $g$ is a $C^\infty$ function on $[a,b]$ and $ g\in C_M([a,b])$. By condition $\mathcal{B}) $, we deduce that the function g is zero and therefore the function f vanishes on the interval $[c, b]$ and so f is identically zero. If $c= a $ [resp. $c= b$]  we set the function
$g = 0$ on $[ a- \eta, a]$ [ resp. $[ b, b+\eta]$ ], where $\eta >0$ and $g = f$ on $[a,b]$.
\end{pre}

Let's eliminate some trivially case of the sequence $M = (M_n)_{n\in\N}$ where the class $ C_M([a,b])$ is quasi-analytic.
\begin{prop}
Let $M=(M_n)_{n\in\N}$ be a sequence of positive numbers. Then the class $ C_M([a,b])$ is quasi-analytic if $\liminf\limits_{n\to\infty}\sqrt[n]{M_n} <\infty$.
\end{prop}
\begin{pre}
Put $\lambda = \liminf\limits_{n\to\infty}\sqrt[n]{M_n}$. For every $\epsilon >0$ corresponds an infinite increasing sequence of natural numbers $(n_j)$ such that $M_{n_j} \leq (\lambda + \epsilon)^{n_j}$. If $ c\in [a,b]$ is such that $ f^{(n)}(c) = 0,\,\,\forall n\in\N$, then $$ |f(x)| \leq \frac{ |f^{(n_j)}(\theta)|}{n_j!}|x- c|^{n_j} \leq \frac{ M_{n_j}}{n_j!}|x- c|^{n_j}\leq e^{(\lambda + \epsilon)}|x- c|^{n_j},\,\,\forall j$$ where $\theta$ is in the open interval of ends $x$  and $c$. Suppose  $|x- c| < 1$, we see then  if we let $j$ tend to infinity, this would imply  $ f=0$ on a subinterval  of $[a,b]$ containing the point $c$. By doing the same thing with one end of this interval and so on, we show that $f$ is zero on $[a,b]$.
\end{pre}
From now on, we assume that
\begin{equation}
\liminf\limits_{n\to\infty}\sqrt[n]{M_n} = \infty.
\end{equation}
Carleman provides a complete answer to the Theorem 2.1.
\begin{thm}
The class  $ C_M([a,b])$ is quasi-analytic if, and only if, $\sum_{n=0}^\infty\frac{1}{\beta_n} =\infty$, where
 $ \beta_n = \inf\limits_{k\geq n}\sqrt[k]{M_k}$.
\end{thm}
For the proof of theorem 2.6, it is often convenient to deal with other equivalent statements involving some other series. In order to do that, we need to change the sequence $M=(M_n)_{n\in\N}$ by an other with suitable properties, and this change does affect   the quasi-analyticity of  the class $ C_M([a,b])$.\\
We introduce here the so-called {\it{ convex regularization by means of the logarithm}}.
\begin{defn}
A sequence of positive real numbers $M=(M_n)_{n\in\N}$ is said to be $log-convex$ if and only if for all $ n \geq 1$ we have that $ M_n^2 \leq M_{n-1}.M_{n+1}$.
\end{defn}
The condition $\liminf\limits_{n\to\infty}\sqrt[n]{M_n} = \infty$ implies the existence of the  convex regularization by means of the logarithm of the sequence $M=(M_n)_{n\in\N}$ see \cite{Man}, that is, a sequence  $M^c = (M_n^c)_{n\in\N}$ such that:
\begin{enumerate}
\item[ a)] the function $k\mapsto \ln M_k^c$ is convex i.e. $ \ln M_k^c \leq \frac{1}{2}( \ln M_{k-1}^c + \ln M_{k+1}^c)$
\item [ b)] $M_k^c \leq M_k,\,\,\,\,\forall k$,
\item [ c)] there is a sequence $ 0 = n_0 < n_1 < n_2 < \ldots $, called the principal sequence, such that $ M_{n_j}^c = M_{n_j}$,
and the function $k\mapsto \ln M_k^c$ is linear in each $[n_j,n_{j+1}]$,

\end{enumerate}
The convex regularized sequence by means of the logarithm is the largest convex minorant of the function $ n\mapsto \log M_n$.
We give an idea of the construction of such sequence.

Let's first recall  the definition of Newton's polygon attached to the sequence $(\log M_n)_{n\in\N}$ (under the condition $\liminf\limits_{n\to\infty}\sqrt[n]{M_n} = \infty)$.\\
Consider in the plane $x0y$ the points $P_n = (n, \log M_n)$. Let   the half-line passing through the point $P_0 = (0, \log M_0)$  and pointing to the negative direction of 0y. Let's turn our half-line in  sense  counter clock wise until it meets a point $P_n = (n, \log M_n)$. Call this point $P_{n_1} = (n_1, \log M_{n_1})$.
The interval $[ P_0,P_1]$  will form the first side of Newton's polygon.
Let's then turn the half-line in the same sense around the point $P_{n_1} = (n_1, \log M_{n_1})$ until it meets a point $P_{n_2} = (n_2, \log M_{n_2})$.
The interval $[ P_1,P_2]$  will form the second side of Newton's polygon, and so on. We thus obtain the Newton's polygon of the sequence
$$ \log M_0, \log M_1,\ldots, \log M_n,\ldots.$$
For all $n\in\N$, let us denote by $(\log M_n)'$ the ordinate corresponding to the abscissa $n$ of the Newton's polygon of the sequence $(\log M_n)_{n\in\N}$.
The sequence $\{(\log M_n)'\}_n$ is the largest convex sequence whose terms are less than the terms of the sequence $(\log M_n)_n$. We have
$$ (\log M_n)' = \inf\limits_{k\geq 0, 0\leq l \leq n}\left(\frac{ k\log M_{n-l} + l\log M_{n+k}}{k+l}\right).$$
Put $M_n^c = \exp (\log M_n)'$, $ M^c = ( M_n^c)_n$, obviously $ M_n^c  \leq  M_n$, $\forall n\in\N$, and
$$ M_n^c = \inf\limits_{k\geq 0, 0\leq l \leq n}\left(  M_{n-l}^{\frac{k}{k+l}} .  M_{n+k}^{\frac{l}{k+l}}\right).$$
It is clear that $ C_{M^c}([a,b]) \subset C_{M}([a,b])$. 
We give some useful properties of $log-convex$ sequences.
\begin{prop}
Assume that $M=(M_n)_{n\in\N}$, with $M_0=1$, is a $log-convex$ sequence of positive real numbers, then
\begin{enumerate}
\item  [a)] the sequence $(\frac{M_{n+1}}{M_n})_{n\in\N}$ is monotone increasing,
\item  [b)] the sequence $(\sqrt[n]{M_n})_{n\in\N}$ is monotone increasing.
\end{enumerate}
\end{prop}
\begin{pre}
\begin{enumerate}
\item [a)]  Is clear from the definition.
\item [b)] $$M_n = \frac{M_n}{M_0} = \prod\limits_{j=1}^n \frac{M_j}{M_{j-1}} \leq (\frac{M_n}{M_{n-1}})^n$$ which gives $M_{n-1}^n \leq M_n^{n-1}$, or equivalently
$\sqrt[n-1]{M_{n-1}}\leq \sqrt[n]{M_n}$.
\end{enumerate}
\end{pre}
Using the  convex regularization by means of logarithm, we give other conditions equivalent to Carleman's one, see \cite{Man}.
\begin{thm}[ Mandelbrojt]
Let $M=(M_n)_{n\in\N}$ be a sequence of positive real numbers. Then the following conditions are equivalent
\begin{enumerate}
\item[i)] $\sum_{n=0}^\infty\frac{1}{\beta_n} =\infty$, where
 $ \beta_n = \inf\limits_{k\geq n}\sqrt[k]{M_k}$
\item[ii)] $\sum_{n=0}^\infty\frac{1}{\sqrt[n]{M^c_n}} =\infty$
\item[iii)] $\sum_{n=0}^\infty\frac{M^c_{n-1}}{M^c_n} =\infty$
\end{enumerate}
\end{thm}
\begin{pre}
Since the sequence $(\sqrt[n]{M^c_n})_{n\in\N}$ is increasing, we have
$$ \sqrt[n]{M^c_n} = \inf\limits_{k\geq n}\sqrt[k]{M^c_k} \leq \inf\limits_{k\geq n}\sqrt[k]{M_k}=\beta_n,$$
hence
\begin{equation}
 \sum\limits_{n=0}^\infty \frac{ 1}{\sqrt[n]{M^c_n}} \geq  \sum\limits_{n=0}^\infty \frac{ 1}{\beta_n}.
\end{equation}
Consider the sequence $(N_n)_{n\in\N}$ defined by
$$ N_{n_i} = M_{n_i}\,\,\,\, \mbox{and}\,\,\,\,N_p = \infty \,\,\,\, \mbox{if}\,\,\,\,p\neq n_i,\,\,\forall i,$$
where $n_1, n_2,\ldots, n_i,\ldots,$ is the principal sequence of $M = (M_n)_{n\in\N}$. Put
$$ \gamma_n = \inf\limits_{ k\geq n} \sqrt[k]{N_k}.$$
We see that $ \gamma_n \geq \beta_n,\,\,\forall n\in\N$ and  if  $ n_{i-1} < n\leq n_i$, we have $\gamma_n = \sqrt[n_i]{M_{n_i}}.$
we then deduce
$$\sum\limits_{n= n_{i-1} + 1}^{n_i} \frac{1}{\gamma_n} = \frac{n_i - n_{i-1}}{ \sqrt[n_i]{M_{n_i}}}.$$
Taking into account that the sequence $(\frac{M^c_{n+1}}{M^c_n})_{n\in\N}$ is monotone increasing, and for all $n_i,\,\,\, M^c_{n_i} = M_{n_i}$, we have
$M_{n_i}^{\frac{1}{n_i}} \leq (\frac{M_{n_i}}{M_{n_{i-1}}})^{\frac{1}{n_i - n_{i-1} }}$, hence
$$ \frac{ \log M_{n_i}}{ n_i} \leq \frac{ \log M_{n_i} - \log M_{n_{i-1}}}{n_i - n_{i-1}}.$$
But we know that the restriction of the function $ p\mapsto \log M_p^c$ to the interval $ [n_{i-1}, n_i]$  is linear and
$$ \forall n\in [n_{i-1}, n_i],\, \,\, \log M_n^c = \frac{ \log M_{n_i} - \log M_{n_{i-1}}}{n_i - n_{i-1}}n + K,$$
for some positive constant $K$. We then obtain
$$\frac{ \log M_{n_i}}{ n_i} \leq \frac{ \log M_{n_i} - \log M_{n_{i-1}}}{n_i - n_{i-1}} = \log M_n^c - \log M_{n-1}^c, \,\, \forall n,\,\,\, n_{i-1} < n \leq n_i.$$
Hence we have: $$ \sqrt[n_i]{M_{n_i}} \leq \frac{M_n^c}{M_{n-1}^c},\,\,\,\forall n,\,\, n_{i-1} < n \leq n_i,$$
 which gives
 $$\sum\limits_{n=1}^\infty \frac{1}{\gamma_n} \geq \sum\limits_{n=1}^\infty \frac{M^c_{n-1}}{M^c_n}.$$
 Since $ \gamma_n \geq \beta_n,\,\,\forall n\in\N$, we have
 \begin{equation}
  \sum\limits_{n=1}^\infty \frac{ 1}{\beta_n} \geq \sum\limits_{n=1}^\infty \frac{M^c_{n-1}}{M^c_n}.
 \end{equation}
 For the last step of the proof we use Carleman's inequality, see \cite{Car} page 112 and \cite{JPW} for the continuous version, that is
 $$ a_1 + \sqrt{a_1a_2}+ \sqrt[3]{a_1a_2a_3}+\ldots + \sqrt[n]{a_1a_2a_3\ldots a_n}+\ldots \leq e (a_1 + a_2+\ldots +a_n+\ldots),$$
 where $a_1,\ldots, a_n,\ldots$ are positives real numbers.\\
 We then see that, if we put $ a_n = \frac{M^c_{n-1}}{M^c_n}$, we obtain
 \begin{equation}
 \sum\limits_{n=1}^\infty \frac{ 1}{\sqrt[n]{M^c_n}}  \leq e\sum\limits_{n=1}^\infty \frac{M^c_{n-1}}{M^c_n}.
 \end{equation}
 Hence the proof the result  by (2.6), (2.7) and (2.8).
\end{pre}
\section{  Proof of sufficiency of theorem 2.6}
We recall some machinery using the theory of metric space, see \cite{Bang}, to prove the sufficiency part of Carleman's theorem (Theorem 2.6).

Let $\mathcal{S}(\R)$ denote the set of all real sequences and $P\subset \N$ an infinite set. We will construct a metric on $\mathcal{S}(\R)$.
\begin{defn}
For any $ X = (x_n)_n\in \mathcal{S}(\R)$ we define
$$ \|X\| = \inf\limits_{k\in P} \left( \max\left( e^{-k} , \max\limits_{ 0\leq n\leq k}|x_n|\right)\right).$$
\end{defn}
If $ X = (x_n)_n\in \mathcal{S}(\R)$, we can assume that $x_0\neq 0$ (if not we change the numbering), let $$\Delta=\{ k\in P\,/\, e^{-k} < |x_0|\},$$
note that $\Delta \neq \emptyset$, hence there exists $ k'\in P$ such that $\inf \Delta =k'$.
\begin{lem}
If $ X = (x_n)_n\in \mathcal{S}(\R)$, then there exists $ 0\leq k\leq k' = \inf \Delta$, $k\in P$, such that
$$ \|X\| = \max \left( e^{-k} , \max\limits_{ 0\leq n\leq k}|x_n|\right).$$
\end{lem}
\begin{pre}
We observe that the sequence $(\max\limits_{0\leq n\leq k}|x_n|)_k$ is increasing in $k$ with infimum $|x_0|$. As the sequence $(e^{-m})_m$ is decreasing, it follows that for all $k\in P$, $ k\geq k'$,
$$ \max \left( e^{-k} , \max\limits_{ 0\leq n\leq k}|x_n|\right) = \max\limits_{0\leq n\leq k}|x_n|,$$
consequently,
$$ \|X\| = \inf\limits_{k\in P,\,\,k\leq k'}\max \left( e^{-k} , \max\limits_{ 0\leq n\leq k}|x_n|\right).$$
The result now follows by  observing that the infimum has been taken over a finite set.
\end{pre}
\begin{lem}
If $ X = (x_n)_n\in \mathcal{S}(\R)$ satisfies the inequality
$$ e^{-k_1} \leq \|X\|\leq e^{-k_2},$$
with $k_1,k_2\in P$. Then there exists $k\in P$, $ k_2 \leq k\leq k_1$ such that $\|X\| = \max \left( e^{-k} , \max\limits_{ 0\leq n\leq k}|x_n|\right)$.
\end{lem}
\begin{pre}
By lemma 3.2, there exists $k\in P$, such that
$ \|X\| = \max \left( e^{-k} , \max\limits_{ 0\leq n\leq k}|x_n|\right)$. Since $\|X\|\leq e^{-k_2}$, we have $e^{-k} \leq e^{-k_2}$, hence $ k\geq k_2$.\\
If $e^{-k_1} \leq \|X\|$, then for all $ p > k_1$, we have $$ e^{-p}  < e^{-k_1} \leq \|X\|= \inf\limits_{k\in P}\left( \max\left( e^{-k} , \max\limits_{ 0\leq n\leq k}|x_n|\right)\right),$$ we obtain then, if $ p > k_1$, $\max\left( e^{-p} , \max\limits_{ 0\leq n\leq p}|x_n|\right) = \max\limits_{ 0\leq n\leq p}|x_n|$. It follows
Then that $$\|X\|= \inf\limits_{k\in P,\, k\leq k_1}\left( \max\left( e^{-k} , \max\limits_{ 0\leq n\leq k}|x_n|\right)\right),$$ which ends  the proof.
\end{pre}
We see that $\|X\| = \|- X\|$, $\|X\| \geq 0$ and $\|X\| = 0$ if and only if $ X =0$ (the zero sequence).
Indeed, it is clear that if $X=0$, then $\|X\|=0$. To prove the converse, suppose that $ X\neq 0$, then there exists $ p\in P$ such that $ |x_p| >0$ and for all $ q\in P$, $ q < p$, $x_q=0$. This implies that
$$\max \left( e^{-j} , \max\limits_{ 0\leq n\leq j}|x_n|\right) \geq |x_p|,\,\,\forall j\in P,\,j\geq p,\,\,\,\mbox{and} \,\,\,\max \left( e^{-j}, \max\limits_{ 0\leq n\leq j}|x_n|\right) = \,\, e^{-j}, \forall j\in P,\,j < p.$$
It then follows that  $\|X\| \geq \min \left( |x_p|,  e^{-p+1}\right) > 0$.\\
We observe that $\|X + Y\| \leq \|X\| + \|Y\|$, for all $ X, Y\in \mathcal{S}(\R)$. If $X =0$ or $Y =0$, the inequality is trivial. Suppose that $ X \neq 0$ and $ Y\neq 0$, then
$$ e^{-k_1} \leq \|X\| \leq e^{-k_2}\,\,\,\,\mbox{and}\,\,\,\,e^{-s_1} \leq \|Y\| \leq e^{-s_2},\,\,\mbox{ for some}\,\,\,\,k_1, k_2,s_1,s_2\in P.$$
By lemma 3.4, there exist  $r, q \in P$, $ k_2\leq r\leq k_1$, $ s_2\leq q\leq s_1$, such that
$$ \|X\|= \max\left( e^{-r} , \max\limits_{ 0\leq n\leq r}|x_n|\right),\,\,\,\,\|Y\|= \max\left( e^{-q} , \max\limits_{ 0\leq n\leq q}|y_n|\right).$$
We can suppose that $ r\leq q$. We have
$$\forall n \leq r,\,\, |x_n + y_n|\leq |x_n| + |y_n|\leq \max\limits_{ 0\leq n\leq r}|x_n| + \max\limits_{ 0\leq n\leq q}|y_n|\leq \|X\| + \|Y\|,$$
hence $ \max\limits_{ 0\leq n\leq r}|x_n + y_n|\leq |X\| + \|Y\|$, we have also $ e^{-r} \leq \|X\| + \|Y\|$, which gives
$$  \|X + Y\| \leq \max\left( e^{-r} , \max\limits_{ 0\leq n\leq r}|x_n + y_n|\right) \leq \|X\| + \|Y\|.$$
We can then provide the space $ \mathcal{S}(\R)$ with a distance function, $d$, defined as follows
$$ d(X, Y)= \|X-Y\|.$$
Let $ X(t)= (x_n(t))_n$, where $ [a,b]\ni t\mapsto x_n(t)$ is a continuous function for each $n\in\N$, we see that the function $t\mapsto \|X(t)\|$ is also a continuous function.

Having developed this machinery, we now utilize it in the proof of the statement that condition iii) of theorem 2.11 implies quasi-analyticity of the relevant functions.\\
Let $M=(M_n)_{n\in\N}$ be a sequence of positive numbers with $M_0 =1$, and suppose that  $\liminf\limits_{n\to\infty}\sqrt[n]{M_n} = \infty$. Consider $M^c = (M_n^c)_{n\in\N}$ the  convex regularization by means of the logarithm of the sequence $M=(M_n)_{n\in\N}$. We denote by $P\subset \N$ the set of all $p\in\N$, such that $ M_p = M^c_p$, we suppose that $0\in P$.
Let  $f\in  C_M([a,b])$ and for each $ t\in [a,b]$, we put $$ X_f(t)=(x_{f,n}(t))_n,\,\,\, \mbox{ where } \,\,\,x_{f,n}(t) = \frac{f^{(n)}(t)}{M_n^ce^{n}},\,\, n\in\N.$$
and
$$ \|X_f(t)\| = \inf\limits_{k\in P}\left( \max\left( e^{-k} , \max\limits_{ 0\leq n\leq k}\frac{|f^{(n)}(t)|}{M_n^ce^{n}}\right)\right).$$
By lemma 3.2, there exists $ l\leq k'$, such that
$$\|X_f(t)\| = \max\left( e^{-l} , \max\limits_{ 0\leq n\leq l}\frac{|f^{(n)}(t)|}{M_n^ce^{n}}\right),$$
where $k'$ is the smallest $p\in P$ such that $ e^{-p} < |f(t)|$. We remark that
$$ e^{-l} \leq \|X_f(t)\|,\,\,\,\,\mbox{and}\,\,\,\, \frac{|f^{(n)}(t)|}{M_n^ce^{n}} \leq \|X_f(t)\|,\,\,\, \forall n= 0, 1,\ldots,l.$$
Suppose that $ t\in [a,b]$ and let $\tau \in\R$ such that $ t +\tau\in [a,b]$. With the index $l$ at the point $ t + \tau$, we have
\begin{equation}
\|X_f(t + \tau)\| \leq \max\left( e^{-l} , \max\limits_{ 0\leq n\leq l}\frac{|f^{(n)}(t + \tau)|}{M_n^ce^{n}}\right).
\end{equation}
The following lemma  gives us a link between $\|X_f(t + \tau)\|$ and $\|X_f(t)\|$.
\begin{lem}[ ]
Suppose that  $\|X_f(t)\| \neq 0$, then
$$ \|X_f(t + \tau)\| \leq \|X_f(t)\| \exp\left( e |\tau|\frac{M^c_l}{M_{l-1}^c}\right),$$
where $l$ satisfies $\|X_f(t)\| = \max\left( e^{-l} , \max\limits_{ 0\leq n\leq l}\frac{|f^{(n)}(t)|}{M_n^ce^{n}}\right)$.
\end{lem}
\begin{pre}
According to above, there  exists $l\in\{1,2,\ldots, k'\}\cap P$ such that
$$\|X_f(t)\| = \max\left( e^{-l} , \max\limits_{ 0\leq n\leq l}\frac{|f^{(n)}(t)|}{M_n^ce^{n}}\right),$$where $k'$ is the smallest $j\in P$ such that $ e^{-j} < |f(t)|$.\\
If $e^{-l} \geq  \max\limits_{ 0\leq n\leq l}\frac{|f^{(n)}(t + \tau)|}{M_n^ce^{n}}$, we have, by (3.1),
$$ \|X_f(t + \tau)\| \leq e^{-l} \leq \|X_f(t)\|,$$
and hence the statement of the lemma holds true.\\
Suppose not, then there exists $0\leq n\leq l$ such that $ \|X_f(t + \tau)\| \leq \frac{|f^{(n)}(t + \tau)|}{M_n^ce^{n}}$.
By using Taylor's theorem for the function $ f^{(n)}$ at the point $t$, we get
\begin{align}
		\begin{split}
			|f^{(n)}(t + \tau)| &\leq \sum\limits_{j=0}^{l-n-1}\frac{|\tau|^j}{j!}|f^{(n+j}(t)| + \frac{|\tau|^{l-n}}{(l-n)!}|f^{(l)}(\xi)|
			  \\
			&\leq \sum\limits_{j=0}^{l-n-1}\frac{|\tau|^j}{j!}M_{j+n}^c e^{n+ j} \|X_f(t)\| + \frac{|\tau|^{l-n}}{(l-n)!} M_l^c e^l\|X_f(t)\|\\
& = \|X_f(t)\| \sum\limits_{j=0}^{l-n}\frac{|\tau|^j}{j!}M_{j+n}^c e^{n+ j}.  \\
		\end{split}
	\end{align}
Then
\begin{align}
		\begin{split}
		\frac{	|f^{(n)}(t + \tau)| }{M_n^c e^n}&\leq \|X_f(t)\|  \sum\limits_{j=0}^{l-n}\frac{|\tau|^j}{j!}(\frac{M_{n+j}^c}{M_{n}^c})  e^j\\
&\leq \|X_f(t)\|  \sum\limits_{j=0}^{l-n}\frac{|\tau|^j}{j!}\left(\frac{M_{l}^c}{M_{l-1}^c}\right)^j  e^j\\
&\leq \|X_f(t)\| \exp\left( e |\tau|\frac{M^c_l}{M_{l-1}^c}\right),
\end{split}
\end{align}
where $log-convexity$ of the sequence $(M_n^c)_n$ has been used to derive $ \frac{M_{n+j}^c}{M_{n}^c}\leq \left(\frac{M_{l}^c}{M_{l-1}^c}\right)^j $.
Hence the lemma is proved.
\end{pre}
Finally we begin the proof of sufficiency of theorem 2.6.\\
Let $ f\in C_M([a,b])$ where $M = (M_n)_n$ with $\liminf\limits_{n\to\infty}\sqrt[n]{M_n} = \infty$. Suppose that $ \sum\limits_{n=0}^\infty \frac{M_n^c}{M_{n+1}^c} =\infty$, where $(M_n^c)_n$ is the  convex regularization by means of logarithm of the sequence $M = (M_n)_n$. Let us always denote by $P \subset\N$  the set of integers $n$ such that $M_n = M_n^c$.
 If in addition there exists $ t_0\in [a,b]$ such that $ f^{(n)}(t_0) =0,\,\,\,\forall n\in\N$, we want to show that $ f \equiv 0$.\\
Let us assume, on contrary, that $ f\neq0$. There exists $ c\in [a,b]$ such that $ 0 = \|X_f(t_0)\| < \|X_f(c)\|$.
Then, there exists $ p\in P$ such that $ 0 = \|X_f(t_0)\| < e^{-p} \leq \|X_f(c)\|$. By the continuity of the function $ t\mapsto \|X_f(t)\|$, there exists $x_0$ between $t_0$ and $ c$ such that $ \|X_f(x_0)\| = e^{-p}$. Successive application of the intermediate value theorem gives a monotonic decreasing sequence of points $ x_0 > x_1 >\ldots, > x_n >\ldots $ such that
$$ \|X_f(x_i)\| = e^{-(p+i)},\,\,\,\,\forall i \geq 0.$$
Sitting $ t = x_i$ and $ t+\tau = x_{i-1}$ and applying lemma 3.6, we get
\begin{equation}
1 \leq e |x_i - x_{i-1}|\frac{M_{l_i}^c}{M_{l_{i-1}}^c},\,\,\,\forall i \geq 1,
\end{equation}
where $l_i$ satisfies
$$ e^{-(p+i)} = \|X_f(x_i)\| = \max \left(e^{-l_i}, \max\limits_{ 0\leq n\leq l_i}\frac{|f^{(n)}(t)|}{M_n^ce^{n}}\right).$$
We note that
$$ e^{-(p+i)} = \|X_f(x_i)\|\leq e^{-(p+i-1)}.$$
By lemma 3.4, we see that $ p+i-1 \leq l_i\leq p+i$. By log-convexity of the sequence $ (M_n^c)_n$, we have $\frac{M_{l_i}^c}{M_{l_{i-1}}^c} \leq \frac{M_{p+i}^c}{M_{{p+i-1}}^c}$.\\
The equation (3.4)can then be rewritten as
$$ \frac{M_{p+i-1}^c}{M_{{p+i}}^c} \leq e |x_i - x_{i-1}|,\,\,\,\,\forall i\geq 1.$$
Summing over all $ i\geq 1$ and using $\sum\limits_{i=1}^\infty |x_i - x_{i-1}|\leq |t_0 - x_0|$, we have
$$ \sum\limits_{i=1}^\infty \frac{M_{p+i-1}^c}{M_{{p+i}}^c} \leq e |t_0 - x_0|.$$
This contradicts the divergence of the series $ \sum\limits_{n=0}^\infty \frac{M_n^c}{M_{n+1}^c}$. Hence it has been proved that $ f\equiv 0$.
\section{Monotonicity property for quasi-analytic Denjoy-Carleman classes}
We know that an analytic function $f$ on the interval $ [a,b]$ is entirely determined by the element $\{f^{(n)}(c)\}_{n\in\N}$, where $ c\in [a,b]$. We are interested in a generalization of this fact, which can be stated as follows:\\
Let $(x_n)_n$  be a sequence of elements of the interval $ [a,b]$. Which condition must check the sequence $(x_n)_n$ in order that the  element $\{f^{(n)}(x_n)\}_{n\in\N}$  determines the function f completely. We see that if the sequence $(x_n)_n$  is constant we get our first property.\\
For an analytic function a response is given by W. Gontcharoff \cite{Gon}.
\begin{thm}\cite{Gon}
Let $f$ be an analytic function on the interval $ [a,b]$. The function $f$ is entirely determined by the knowledge of the values  $f^{(n)}(x_n),\,\,\,(n=1,2\ldots)$ if the series  $\sum\limits_{n=1}^\infty |x_{n-1} - x_n|$ converges.
\end{thm}
As a consequence, we deduce that if $f^{(n)}(x_n) = 0,\,\,\,\forall n\in\N$,  and the series $\sum\limits_{n=1}^\infty |x_{n-1} - x_n|$ converges, then the function $f$ is identically zero.\\

The question now is whether a similar result remains valid for a quasi-analytic class. In the case of a quasi-analytic class of Denjoy-Carleman, we have the following theorem proved by W.Bang \cite{Bang}.\\
Let $M = (M_n)_n$ be a sequence of positive numbers with $M_0= 1$. Suppose that $M = (M_n)_n$ is logarithmically convex and satisfying one of the equivalent conditions of theorem 2.11.

\begin{thm}
Let $f\in C^\infty([a,b])$ satisfy
$$ \sup\limits_{t\in [a,b]}|f^{(n)}(t)| \leq M_n,\,\,\, n\in\N.$$
Suppose that there exists a sequence $(x_n)_n$, $ x_n\in[a,b]$, such that $f^{(n)}(x_n) = 0,\,\,\,\forall n\in \N$. If the series $\sum\limits_{n=1}^\infty |x_{n-1} - x_n|$ converges, then $f$ is identically $0$.
\end{thm}

For the convenience of the reader and for completeness, we reproduce the proof.

We consider, for each $n\in\N$, the function
$$ B_{f,n}(t) = \sup\limits_{ j\geq n} \frac{| f^{(j)}(t)|}{e^j M_j},\,\,\,\, t\in [a,b].$$
We list some properties of this sequence of functions.
\begin{enumerate}
\item $B_{f,n}(t)\leq e^{-n},\,\,\,\,\forall t\in [a,b]$,
\item $ B_{f,0}(t)\geq B_{f,1}(t)\geq \ldots \geq B_{f,n}(t)\geq \ldots ,$
\item if $f^{(n)}(t_0) =0$, then $ B_{f,n}(t_0) = B_{f,n+1}(t_0)$.
\end{enumerate}
\begin{lem}
The function $ t\mapsto B_{f,n}(t)$ satisfies the estimate
$$ B_{f,n}(t + \tau)\leq \max\left( B_{f,n}(t), e^{-q}\right)\exp\left( e|\tau|\frac{M_q}{M_{q-1}}\right),$$
for every $ q\in\N,\,\,\,\, q > n$ and $ t, t+\tau \in [a,b]$.
\end{lem}
\begin{pre}
We follow the proof of the lemma 3.6. Let $ j \in\N,\,\,\,\, n \leq j < q$
\begin{align}
		\begin{split}
		\frac{|f^{(j)}(t + \tau)|}{e^j M_j} &\leq \sum\limits_{i=0}^{q-j-1}\frac{|\tau|^i}{e^j M_j i!}|f^{(j+i)}(t)| + \frac{|f^{(q)}(\xi)||\tau|^{q-j}}{e^jM_j(q-j)!}\\
&= \sum\limits_{i=0}^{q-j-1}\frac{ M_{i+j}}{M_j}\frac{|f^{(j+i)}(t)|}{M_{i+j}e^{i+j}}\frac{(e|\tau|)^i}{ i!} +
e^{-q}\frac{M_q}{M_j}\frac{|f^{(q)}(\xi)|}{M_q}\frac{(e|\tau|)^{q-j}}{(q-j)!}\\
& \leq B_{f,n}(t)\sum\limits_{i=0}^{q-j-1}\frac{\left(e|\tau|\right)^i}{ i!}(\frac{ M_q}{M_{q-1}})^i + e^{-q}\left(\frac{ M_q}{M_{q-1}}\right)^{q-j}
\frac{(e|\tau|)^{q-j}}{ (q-j)!}\\
& \leq \max\left( B_{f,n}(t), e^{-q}\right)\exp\left( e|\tau|\frac{M_q}{M_{q-1}}\right).
\end{split}
\end{align}
Where we used that the sequence $ M=(M_n)_n$ is logarithmically convex.\\
As a consequence we see that the function $ t\mapsto B_{f,n}(t)$ is continuous on the interval $[a,b]$. Indeed, if $t_0\in [a,b]$, there exists $ q\in\N$, such that $ e^{-q} < B_{f,n}(t_0)$. By lemma 4.3, we have
$$  \big| B_{f,n}(t_0 + \tau) - B_{f,n}(t_0)\big|\leq B_{f,n}(t_0)\left( \exp\left(  e|\tau|\frac{M_q}{M_{q-1}}\right)-1\right).$$
Let $f$ and $(x_j)_j$ be as in the theorem. Set $ \tau_k = \sum\limits_{j=0}^{k - 1} |x_j - x_{j+1}|,\,\,\,\, k\geq 1$ and $\tau_0 = 0$.
We remark that if $ t\in [ \tau_{n-1}, \tau_n]$, then
$$
\left\lbrace
\begin{array}{ll}
x_{n - 1} +\tau_{n-1} -t\in [ x_n, x_{n-1}]\subset [ a, b] & \mbox{if $ x_n < x_{n-1}$ }\\
x_{n - 1}  -\tau_{n-1} + t\in [x_{n-1}, x_n ]\subset [ a, b] & \mbox{if $ x_n  > x_{n-1}$ }
\end{array}
\right.
$$
We define a function $\overline{B}_{f,n}$ on $[ \tau_{n-1}, \tau_n]$ by
$$ \overline{B}_{f,n}(t) = \left\lbrace
\begin{array}{ll}
B_{f,n}(x_{n - 1} +\tau_{n-1} -t)  & \mbox{if $ x_n < x_{n-1}$ }\\
B_{f,n}(x_{n - 1} -\tau_{n-1} + t) & \mbox{if $ x_n  > x_{n-1}$ }
\end{array}
\right. $$
The function $ t\mapsto \overline{B}_{f,n}(t)$ is continuous and by 3)
\begin{equation}
\overline{B}_{f,n}(\tau_n) = B_{f,n}(x_n) = B_{f,n+1}(x_n)=\overline{B}_{f,n+1}(\tau_n).
\end{equation}
We can then paste together the functions  $ [ \tau_{n-1}, \tau_n] \ni t\mapsto \overline{B}_{f,n}(t)$ with different $n\in\N$ and define a new  function
on the interval $[0, \tau[$, where $ \tau = \sup\limits_{n}\tau_n$, by
$$  b_f(t)=\overline{B}_{f,n}(t),\,\,\,\,\mbox{ if} \,\,\,\, t\in [ \tau_{n-1}, \tau_n].$$
This is a continuous function. By 1. and 2. we find that
\begin{equation}
 b_f(t) \leq e^{-n},\,\,\,\, \forall t \geq \tau_{n-1},
\end{equation}
 and hence $b_f(t)\rightarrow 0$ as $t\rightarrow \tau$. Note that for all $ t_1, t_2\in [ \tau_{n-1}, \tau_n]$, we have
\begin{equation}
\forall q\geq n,\,\,\,\,\,b_f(t_2)\leq \max\left(b_f(t_1), e^{-q}\right) \exp\left(  e|t_1-t_2|\frac{M_q}{M_{q-1}}\right).
\end{equation}
Indeed, since $ t_1, t_2\in [ \tau_{n-1}, \tau_n]$,  we have
$$ b_f(t_2) = \overline{B}_{f,n}(t_2) = \left\lbrace
\begin{array}{ll}
B_{f,n}(x_{n - 1} +\tau_{n-1} -t_2)  & \mbox{if $ x_n < x_{n-1}$ }\\
B_{f,n}(x_{n - 1} -\tau_{n-1} + t_2) & \mbox{if $ x_n  > x_{n-1}$ }.
\end{array}
\right. $$
Using the lemma 4.3,  we deduce (4.4).

If $f\neq 0$, then $b_f\neq 0$, hence $]0,\mu]\subset Im\, b_f$, for some $\mu >0$. Let $k_0$ be the smallest integer such that $ e^{-k_0}\in ]0,\mu]\subset Im \, b_f$. Put
$$ t_{k_0} = \inf\{ t\in [0,\tau[/\,\,\, b_f(t)= e^{-k_0}\},\,\,\,\mbox{and for }\,\,\, k >k_0,\,\,\,t_k = \inf\{ t\in ]t_{k-1},\tau[/\,\,\, b_f(t)= e^{-k}\}.$$
So we have a strictly increasing sequence $(t_k)_{k\geq k_0}$ such that $ b_f(t_k)= e^{-k}$. We see that for each $ t\in ]t_{k-1}, t_k[$, $b_f(t) > e^{-k}$, and by (4.3), we have $ t_k < \tau_k,\,\,\,\forall k\geq k_0$.\\
The  sequence $(\tau_n)_{n\in\N}$ defines on each interval $]t_{k-1},t_k[$ a subdivision. More precisely, Let $s\in\N$ be the largest integer such that $\tau_s \leq t_{k-1}$ and $r\in\N$ the smallest integer such that $t_k\leq \tau_{r}$. Using the fact that for each $ t\in ]t_{k-1}, t_k[$, we have $b_f(t) > e^{-k}$ and (4.4), we obtain

\begin{equation}
	\left\lbrace
		\begin{aligned}
			b_f(t_{k-1})	&\leq b_f(\tau_{s+1})\exp\left( e(\tau_{s+1} - t_{k-1})\frac{M_{k}}{M_{k-1}}\right)\\
			b_f(\tau_{s+1})	&\leq b_f(\tau_{s+2})\exp\left( e(\tau_{s+2} - \tau_{s+1})\frac{M_{k}}{M_{k-1}}\right)\\
			\vdots  &\hspace{2cm}\vdots \\
            b_f(\tau_{k-1})	&\leq b_f(t_k)\exp\left( e(t_k - \tau_{k-1}) \frac{M_{k}}{M_{k-1}}\right)
		\end{aligned}
	\right.
\end{equation}
We infer that 
$$b_f(t_{k-1}) \leq b_f(t_k) \exp\left( e(t_k - t_{k-1}) \frac{M_{k}}{M_{k-1}}\right).$$
 Since $b_f(t_k) = e^{-k}$ for every $k \geq k_0$, we have
 $$ t_k - t_{k-1}\geq \frac{1}{e}\frac{M_{k-1}}{M_{k}},\,\,\,\,\,\forall k\geq k_0.$$
Hence 
$$ t_k \geq t_{k_0} + \frac{1}{e}\sum\limits_{j= k_0 + 1}^k\frac{M_{j}}{M_{j-1}}.$$
Since  $ \tau_k > t_k$, we obtain
$$  \sum\limits_{j=0}^{k - 1} |x_j - x_{j+1}|> t_{k_0} + \frac{1}{e}\sum\limits_{j= k_0 + 1}^k\frac{M_{j}}{M_{j-1}}.$$
which proves the theorem.
\end{pre}
As consequence of this theorem, we have the following result:\\
Let $M = (M_n)_n$ be a sequence of positive numbers with $M_0= 1$. Suppose that $M = (M_n)_n$ is logarithmically convex and satisfying one of the equivalent conditions of theorem 2.11.
\begin{cor}
	Let $f\in C^\infty([a,b])$ satisfy
	$$ \sup\limits_{t\in [a,b]}|f^{(n)}(t)| \leq M_n,\,\,\, n\in\N.$$ 
	If $f^{(n)}(a)>0$ for all $n\in\N$, then $f^{(n)}(x)>0$ for all $n\in\N$ and all $x\in[a,b]$.
\end{cor}
\begin{pre}
Suppose that $f^{(n)}(a)>0$ for all $n\in\N$ and there exists $k_0\in\N$ such that $f^{(k_0)}$ has a zero 
$x_{k_0}\in ]a,b]$.
Then there exists $ x_{k_0+1} <  x_{k_0}$ such that 
$f^{(k_0+1)}(x_{k_0+1})=0$. Continuing, we find a strictly decreasing sequence $x_{k_0} > x_{k_0+1}> \ldots$ where $x_l$ is a zero of $f^{(l)}$, for all $ l\geq k_0$. By theorem 4.2, the function $f^{(k_0)}$  is identically $0$, hence $ f$ is the restriction to the interval $[a,b]$ of a polynomial of degree at most $k_0-1$, which is a contradiction, since $f^{(n)}(a)\neq 0$ for all $n\in\N$.\\

As an immediate consequence of this corollary, we deduce, according to Bernstein's theorem, see 
\cite{Widder}, page 146, that the function $f$ can be extended analytically into the plane of complex numbers to a holomorphic function in the disk $\left\lvert z -a\right\rvert < b -a$.
 \end{pre}
 \begin{thm}[Carleman]
Let $ C_M([a,b])$ be a quasianalytic class which contains strictly the analytic class. Then the Borel mapping $$T_c: C_M([a,b])\to \R[[x]],\,\,\,\,\,c\in [a,b],\,\,\,T_c =\sum\limits_{n=0}^\infty \frac{f^{(n)}(c)}{n!}x^n. $$ 
is  not surjective.
\end{thm}
\begin{pre}
	Consider a non-convergent series 
	$\sum\limits_{n\in\N} a_n x^n$ such that $ a_n >0,\,\,\,\,\forall n\in\N$. By corollary 4.5, we see that $\sum\limits_{n\in\N} a_n x^n\notin T_c\left(C_M([a,b])\right)$. 
\end{pre}	
\section{Quasi-analytic classes associated to a sequence of integers  }
  For a function $f$ which is  $C^\infty$   on $[a,b]$, we can formulate the principle of Pringsheim as follows:
\begin{equation}
f\,\,\,\mbox{ is analytic on the interval}\,\,\,\,[a,b]     \Leftrightarrow\,\,\limsup\limits_{n\to \infty}\frac{1}{n}\sqrt[n]{M_n(f)} < \infty\,.
\end{equation}
To preserve the validity of the condition $\mathcal{DC}$) for a function $f$, one only needs a weakened version of the principle of Pringsheim, namely the following condition:
\begin{equation}
\liminf\limits_{n\to \infty}\frac{1}{n}\sqrt[n]{M_n(f)} < \infty\,.
\end{equation}
More precisely we have:
\begin{prop}
Let $f$ be a $ C^\infty$ function on $[a,b]$ such that $$\liminf\limits_{n\rightarrow \infty}\frac{\sqrt[n]{M_n(f)}}{n} < \infty .$$ If $ c \in [0,1]$ is such that $ f^n(c) = 0$, $\forall n\in \N$, then $f$ is identically null on $[a,b]$.
\end{prop}

\begin{pre}
By Taylor's formula, we have, for each $n\in\N$,
\[ f(x)= \frac{f^{(n)}\left(c+ \theta_n(x-c)\right)}{n!} (x-c)^n,\]
where $0 < \theta_n < 1$. Let $(n_k)_k$ be an infinite subsequence such that:
\[ \lim\limits_{ n_k \to \infty}  \frac{\sqrt[n_k]{M_{n_k}(f)}}{n_k}       =        \liminf\limits_{n\rightarrow \infty}\frac{\sqrt[n]{M_n(f)}}{n}.\]

Hence, there exists $ A >0$, such that $$M_{n_k}(f) \leq A^{n_k} \factorial{n_{n_k}}, \,\,\,\,\,\forall  n_k.$$ We have then
$$ f^{(n_k)}\left(c+ \theta_{n_k}(x-c)\right)\leq M_{n_k}(f) \leq A^{n_k}\factorial{ n_{n_k}},$$
and consequently
$$ |f(x)| \leq (A |x-c|)^{n_k},\,\,\,\forall k\in\N.$$
Hence if $|x-c| < \frac{1}{A}$ its follows that $f(x)=0$.\\
Choosing in place of $c$ the point $ c\pm \frac{1}{2A}$ and once more repeating the same reasoning, we obtain $f(x)= 0$ on the whole interval $ [a,b]$.
\end{pre}
Let $\overline{n}= (n_k)_{k\in\N}$ be a strictly increasing sequence of natural numbers, we denote by $C_{\overline{n}}([a,b])$ the set of all $C^\infty$ functions on $[a,b]$ such that, there exist two positive
constants $A,\,B$, such that 
\begin{equation}
 {M_{n_k}(f)} \leq B A^{n_k} \factorial{n_k},\,\,\,\,\forall k\in\N,
 \end{equation}
By proposition 5.1,  any sequence of natural numbers $\overline{n}= (n_k)_{k\in\N}$, which increases without limit, defines some quasi-analytic class of functions which satisfies the condition $\mathcal{DC})$. We will call such class, quasi-analytic class  with respect to the sequence $\overline{n}= (n_k)_{k\in\N}$.
\section{Monotonicity property for quasi-analytic  classes associated to a sequence of integers}
\subsection{Generalized  Taylor's Theorem}
Let $(x_n)_{n\in\N}$ be a sequence of real numbers. For each $n\geq 1$, there is a polynomial $Q(x,x_0,x_1,\ldots,x_{n-1})$ of degree $n$ defined by 
\begin{equation}
\left\lbrace
\begin{array}{l}
Q^{(m)}(x_m,x_0,x_1,\ldots,x_{n-1}) =0,\,\,\, m = 1,\ldots, n-1,    \\
Q^{(n)}(x,x_0,x_1,\ldots,x_{n-1}) =1.
\end{array}
\right. 
\end{equation}
By definition, if $n =0$, we put $Q(x)\equiv 1$.\\
To obtain this polynomial, for $n \geq 1$, one carries out $n $ indefinite integration of the unit function and determining the constants so that conditions (6.1) are verified. We find the expression
\begin{equation}
Q(x,x_0,x_1,\ldots,x_{n-1}) = \int_{x_0}^{x}  \,\mathrm{d}t_1 \int_{x_1}^{t_1}  \,\mathrm{d}t_2 \int_{x_2}^{t_2}  \,\mathrm{d}t_3\ldots \int_{x_{n-1}}^{t_{n-1}}  \,\mathrm{d}t_n.
\end{equation}
Note that by condition (6.2), we have 
\begin{equation}
Q(x,x_0,x_1,\ldots,x_{n-1}) = \int_{x_0}^{x} Q(t,x_1,x_2,\ldots,x_{n-1})\mathrm{d}t,
\end{equation}
and 
for all $m\leq n$
$$ Q^{(m)}(x,x_0,x_1,\ldots,x_{n-1}) = Q(x,x_m,\ldots,x_{n-1}).$$
 We say that the system of polynomials $ Q(x,x_0,x_1,\ldots,x_{n-1}),\,\, n = 0, 1, 2,\ldots $ is associated to the sequence $(x_n)_{n\in\N}$. Using (6.3), polynomials $ Q(x,x_0,x_1,\ldots,x_{n-1}),\,\, n = 0, 1, 2,\ldots $ can be calculated step by step:
 $$
\left\lbrace
\begin{array}{ll}
Q(x,x_0)	&=\left( x-x_0\right)\\
Q(x,x_0,x_1)	&=\frac{1}{2}\left(\left(x - x_1\right)^2 - \left(x_0 - x_1\right)^2\right)\\
Q(x,x_0,x_1,x_{2})	&=\frac{1}{\factorial{3}}\left(\left(x - x_2\right)^3 - 3\left(x_1 - x_2\right)^2	\left( x-x_0\right) -\left(x_0 - x_2\right)^2\right)
\end{array}
\right.
$$
 For any integer $n >0$, the polynomial $ Q(x,x_0,x_1,\ldots,x_{n-1})$ satisfies the following relations 
 \begin{align}
	\begin{split}
Q(x_0,x_0,x_1,\ldots,x_{n-1})&=0,\\
Q'(x_1,x_0,x_1,\ldots,x_{n-1})  & = 0,\\ 
Q^{(n-1)}(x_{n-1},x_0,x_1,\ldots,x_{n-1}) & =0,\\
Q^{(n)}(x_n,x_0,x_1,\ldots,x_{n-1})& = 1.
	\end{split}
\end{align}
 
Suppose now that $x_0 > x_1 >\ldots >x_n >\ldots$.\\
Taking into account that
$ x_n < x_j,\,\,\,j= 1,\ldots, n-1$ and $x_0 > x_1,\ldots > x_{n-1}$, we have 
$$\int_{x_0}^{x}  \,\mathrm{d}t_1 \int_{x_1}^{t_1}  \,\mathrm{d}t_2 \int_{x_2}^{t_2}  \,\mathrm{d}t_3\ldots \int_{x_{n-1}}^{t_{n-1}}  \,\mathrm{d}t_n \leq \int_{x_n}^{x}  \,\mathrm{d}t_1 \int_{x_n}^{t_1}  \,\mathrm{d}t_2 \int_{x_n}^{t_2}  \,\mathrm{d}t_3\ldots \int_{x_{n}}^{t_{n-1}}  \,\mathrm{d}t_n,$$
and 
$$
\int_{x_0}^{x}  \,\mathrm{d}t_1 \int_{x_0}^{t_1}  \,\mathrm{d}t_2 \int_{x_0}^{t_2}  \,\mathrm{d}t_3\ldots \int_{x_{0}}^{t_{n-1}}  \,\mathrm{d}t_n\leq
\int_{x_0}^{x}  \,\mathrm{d}t_1 \int_{x_1}^{t_1}  \,\mathrm{d}t_2 \int_{x_2}^{t_2}  \,\mathrm{d}t_3\ldots \int_{x_{n-1}}^{t_{n-1}}  \,\mathrm{d}t_n 
$$
Which gives, by (6.2),
	\begin{equation}
		\frac{(x-x_0)^n}{\factorial{n}} \leq Q(x,x_0,x_1,\ldots,x_{n-1}) \leq \frac{(x-x_n)^n}{\factorial{n}}.
		\end{equation} 
We see that for all $x > x_0$, we have 	$ Q(x,x_0,x_1,\ldots,x_{n-1}) > 0$.

Let $ f\in C^\infty \left([a,b]\right)$ and suppose that for all $n\in\N$, $x_n\in [a,b]$ with $ x_0 < b$. We put, for each $n\in\N$,
\begin{equation} 
	\begin{split}                                  
		R_n(f)(x)&=\\      
		&f(x)-f(x_0)- f'(x_1)Q(x,x_0)-\dots - f^{(n)}(x_n)Q(x,x_0,\ldots,x_{n-1}). 
	\end{split}                                        
\end{equation} 
For all $ 0\leq k\leq n $, we have, by (6.4),
$$\left(R_n(f)\right)^{(k)}(x_k) = f^{(k)}(x_k) - f^{(k)}(x_k)Q^{(k)}(x_k,x_0,x_1,\ldots,x_{k-1}) =0.
$$
and 
$$ \left(R_n(f)\right)^{(n+1)}(t) = f^{(n+1)}(t),\,\,\,\forall t\in [a,b].
$$
We deduce then
$$ R_n(f)(x) = \int_{x_0}^{x}  \,\mathrm{d}t_1 \int_{x_1}^{t_1}  \,\mathrm{d}t_2 \ldots \int_{x_{n-1}}^{t_{n-1}}  \,\mathrm{d}t_n
\int_{x_{n}}^{t_{n}}  \left(R_n(f)\right)^{(n+1)}(t)\,\mathrm{d}t.
$$
Hence there exists $ \xi\in ]x_n,x[$ such that
\begin{equation}
	R_n(f)(x) = f^{(n+1)}(\xi) Q(x,x_0,x_1,\ldots,x_n).
\end{equation}	
We have shown the following proposition:
\begin{prop}
	Let $ f\in C^\infty \left([a,b]\right)$, then, for each $m\in\N$,
	\begin{equation}
	\begin{split}                                  
	f(x)&=\\      
		&f(x_0) + f'(x_1)Q(x,x_0) +\dots + f^{(m)}(x_m)Q(x,x_0,\ldots,x_{m-1}) + f^{(m+1)}(\xi) Q(x,x_0,x_1,\ldots,x_m), 
	\end{split} 
\end{equation}	
where $ \xi\in ]x_n,x[$.
\end{prop}
Let $\overline{n}= (n_k)_{k\in\N}$ be a strictly increasing sequence of natural numbers, and  $f\in C_{\overline{n}}([a,b])$.
\begin{thm}
Suppose that $ f^{(n)}(x_n)= 0$, for all $ n\in\N $, where the sequence $ (x_n)_{n\in\N}$ is as above.
If the series $\sum\limits_{j=0}^\infty |x_j - x_{j+1}|$ converges, then the function $f$ is identically null.	

\end{thm}	
\begin{pre}
Since 	$f\in C_{\overline{n}}([a,b])$, there exist two positive constants $A,\,B$, such that 
$$
 {M_{n_k}(f)} \leq B A^{n_k} \factorial{n_k},\,\,\,\,\forall k\in\N,
 $$
Put $ R_p = \sum\limits_{j=p}^\infty |x_j - x_{j+1}|$. Let $q_0\in\N$, such that for all $q\geq q_0, \,\,\, R_q < \frac{1}{A}$. We define the sequence $ (m_s)_{s\in\N}$, by 
$$ \forall s\in\N,\,\,\, m_s+q+1\in\{ n_k\,/\,k\in\N\}.$$
It is clear that $m_s\to\infty$ when $ s\to\infty$.

We apply Proposition 6.1 for $ m_s+q$:	
$$f(x) = f^{(m_s+q+1)}(\xi) Q(x,x_0,x_1,\ldots,x_{m_s+q}).$$
Taking $q$ times the derivative of each member, we obtain:
$$f^{(q)}(x) = f^{(m_s+q+1)}(\xi) Q^{(q)}(x,x_0,x_1,\ldots,x_{m_s+q}) = f^{(m_s+q+1)}(\xi) Q(x,x_{q},x_{q+1},\ldots,x_{m_s+q}).$$ 
We deduce
$$ |f^{(q)}(x)|\leq B A^{q+1} \frac{\factorial{(m_s+q+1)}}{\factorial {m_s}}\left(A(|x-x_q|+ R_q)\right)^{m_s},$$  
we have used $|x-x_{m_s+q+1}|\leq |x-x_q|+ R_q$. \\If 
$\mu = \frac{1- AR_q}{A}$, we see that $\forall x\in [x_q -\mu, x_q +\mu],\,\,\, A\left(|x-x_q|+ R_q\right) < 1$ and therefore the function $f^{(q)}$ is zero on the interval $[x_q -\mu, x_q +\mu]$. Hence $f^{(q)}$ is zero on the interval $[a, b]$, by Proposition 5.1. If $ q=1$, the result has been proved. If $ q > 1$, the function $f^{(q-1)}$  is constant and since $f^{(q-1)}(x_{q-1}) = 0$ we deduce that $f^{(q-1)}$ is zero on the interval $[a, b]$. Continuing, we find that $f$ is zero on the interval $[a, b]$.
\end{pre}
\begin{rem}
As a consequence of the Theorem 6.2, we deduce that the corollaries 4.5 and 4.6 are also valid for $C_{\overline{n}}([a,b])$.
\end{rem}
\begin{prob}
	Does any quasianalitic class satisfy the property of monotonicity?
	To fix the ideas, let $\mathcal{R} $ be a polynomially bounded o-minimal structure that extends the structure defined by the restriction to the interval $ [a,b]$ of  analytic functions, see, \cite {VenDen1}. Is the monotonicity property verified  by the definable functions in this structure? We know an answer for the structure defined by the multisommables functions, see \cite {VenDen}, corollary 8.6.
\end{prob}	

	\bf{
 Abdelhafed Elkhadiri\\
University Ibn Tofail,
Faculty of Sciences\\
Kenitra, Morocco\\
E-mail: elkhadiri.abdelhafed@uit.ac.ma}
\end{document}